\documentclass[11pt,a4paper]{article}
\usepackage[utf8]{inputenc}  
\usepackage[]{times}
\usepackage{epsfig}
\usepackage{psfrag}
\usepackage{epsfig}
\usepackage{indentfirst}
\usepackage{graphicx}
\usepackage{harvard}
\usepackage{algorithm}
\usepackage{algorithmic}
\usepackage{multirow}
\usepackage[english]{babel}
\usepackage{amssymb}
\usepackage{amsmath}
\usepackage{bchart}
\usepackage{pgfplots}
\usepackage{hhline}
\usepackage{subfigure}
\usepackage{times}

\newcommand{\citeny}[1]{\citename{#1}~\citeyear{#1}}
\newcommand{\ipre}{\lessdot}

\usepackage{todonotes}

\begin{document}

\bibliographystyle{dcu}
\pagestyle{empty}

\setlength{\itemsep}{0pt}

\title{\textbf{A new look at the bowl phenomenon}}
\author           {Pedro B. Castellucci$^1$ and Alysson M. Costa$^2$\\ \\
\small \emph{$^1$ Instituto de Ciências Matemáticas e de Computação, Universidade de São Paulo} \\
\small \emph{$^2$ Department of Mathematics and Statistics, University of Melbourne} \\
\small pbc@icmc.usp.br, alysson.costa@unimelb.edu.au}
\date{}
\maketitle

\thispagestyle{empty}

\vspace{1cm}

\vspace{0.5cm}

\noindent \textbf{Abstract} \\
An interesting empirical result in the assembly line literature states that slightly unbalanced assembly lines (in the format of a bowl - with central stations less loaded than the external ones) present higher throughputs than perfectly balanced ones. This effect is known as the \emph{bowl phenomenon}. In this study, we analyze the presence of this phenomenon in assembly lines with integer task times. For this purpose, we modify existing models for the \emph{simple assembly line balancing problem} and \emph{assembly line worker assignment and balancing problem} in order to generate configurations exhibiting the desired format. These configurations are implemented in a stochastic simulation model, which is run for a large set of recently introduced instances. The obtained results are analyzed and the findings obtained here indicate, for the first time, the existence of the bowl phenomenon in 
a large set of configurations (corresponding to the wide range of instances tested) and also the possibility of reproducing such phenomenon in lines with a heterogeneous workforce.
{\small }



\newpage



\citationmode{abbr}

\section{Introduction}
\label{sec:introducao}
Assembly lines are productive systems specially useful for large scale standardized manufacturing. The main rationale behind such systems is labour division. This is done by partitioning  the tasks to be executed among a number of workers or workstations \cite{scholl99balancing}. In each workstation, a subset of the tasks is executed. The   classical assembly line balancing problem is known as SALBP (\emph{simple assembly line balancing problem}) and considers, among other hypotheses, that all workers/workstations are equally efficient. A large amount of research has been conducted on the SALBP and includes both classical and recent surveys and research articles \cite{salveson55assembly,tonge61heuristic,baybars86survey,ghosh89comprehensive,scholl99balancing,scholl06state,boysen08versatile,scholl10absalom}.  

SALBP's hypotheses are rarely valid in practical contexts. This fact has motivated the study of a large number of its variants \cite{boysen07classification,boysen08assembly,battaia12taxonomy}. In particular, this research addresses, besides SALBP, the case of assembly lines in sheltered work centers for disabled. In this case, due to workers heterogeneity, task execution times depend on the worker they are assigned to. This gives origin to a problem in which, besides assigning tasks to workstations, there is also the need to assign workers to workstations. This situation has been first modelled by \citeny{miralles07advantages} and named ALWABP (\emph{assembly line worker assignment and balancing problem}).

Both SALBP and ALWABP assume deterministic execution times for each task and, hence, the line's throughput is dictated by the most loaded workstation. In practice, however, task execution times are usually stochastic. This may occur due to small differences among parts, changes in worker's behavior, technical issues and many other random events. The randomness in execution times may lead to interesting effects such as the well known bowl phenomenon, in which slightly unbalanced lines in the format of a bowl (i.e.,  with central workstations less loaded than external ones), or with central workstations with more stable loads (less variable execution times), have  higher productivities \cite{hillier66effects}.

The theory behind the Bowl phenomenon is not yet fully understood. Nevertheless, the effect may lead to adjustments in optimization methods (heuristic and exact) in order to enhance throughput of assembly lines in practical contexts. Thereunto, simulation models are helpful, because they allow the modelling of execution times with higher precision.

In this study, an assembly line simulation model for productivity analysis is proposed. In particular, the interest is focused on recently proposed instances for the SALBP \cite{otto13salpbgen} and marginally on instances of the ALWABP. The objective is twofold:

Our first goal is to verify the possibility of generating solutions for the ALWABP that can take advantage of the bowl phenomenon. For this purpose, a mixed-integer programming model capable of generating solutions with unbalanced bowl-shaped profiles (with greater or lesser slopes) is proposed. These solutions are used as parameters in the simulation model, also developed in this work, and the results are discussed.

The second and more general goal is to evaluate the own existence of the bowl phenomenon. Indeed, the literature has always dealt with a reduced number of configurations (obtained from existing classical case studies). Here, we use the newly proposed instances of \cite{otto13salpbgen} and effect a much larger computational study. Also, since the parameters for the simulation model come from the proposed mixed-integer programming model, the indivisibility of tasks is considered for evaluation which is not usually reported in the literature. 

The remainder of this article is organized as follows. In the following section, SALBP and ALWABP are discussed and their models are presented. Also in this section, these model are modified in order to generate solutions with unbalanced bowl-shaped profiles. Then, Section~\ref{sec:modelo} focus on the proposed simulation model development while Section~\ref{sec:resultados} presents the computational experiments and results. General conclusions end this paper in Section~\ref{sec:conclusoes}.

\section{SALBP and ALWABP}
\label{sec:definicoes}

The most idealized problem of assembly line balancing is the \textit{simple assembly line balancing problem}, SALBP. In this problem, a paced assembly line with a fixed cycle time is considered. Task execution times are deterministic, task allocations are only constrained by precedence relations, lines are serial and all the workstations are equally equipped \cite{scholl99balancing}.

In order to model the SALBP, let $N$ be a set of tasks. Tasks $i\in N$ must respect a partial ordering and the execution time of each task $i\in N$ is $t_i$. We use the notation $i\ipre j$ to indicate that task $i$ must be completed before task $j$ is initiated. Defining binary variables $x_{ik}$ (equal to 1 if, and only if, task $i$ is executed at station $k$). Let $C$ be the cycle time and $S$ the set of workstation, then a SALBP model can be written as:

\begin{equation}
\label{salbp_1}
\textrm{Min }   C
\end{equation}
\hspace{1cm} subjected to:
\begin{alignat}{2}
\label{salbp_2}
\sum_{k \in S}  x_{ik} = 1, &&  \qquad &\forall i \in N,\\
\label{salbp_3_cicle}
\sum_{i\in N} t_{i}\cdot x_{ik} \leq C, && \qquad &\forall k \in S,\\
\label{salbp_4}
\sum_{k\in S}  k\cdot x_{ik} \leq  \sum_{k\in S} k\cdot x_{jk}, && \qquad & \forall i,j \in N, i\ipre j,\\
\label{salbp_5}
x_{ik} \in \{0,1\}, && \quad &\forall i\in N, \forall k\in S.
\end{alignat}

The goal of the model (\ref{salbp_1} - \ref{salbp_5}) is to provide an allocation of tasks to workstation so that the cycle time is minimized. Constraints (\ref{salbp_2}) force every task to be executed, while constraints (\ref{salbp_3_cicle}) define the assembly line's cycle time and (\ref{salbp_4}) ensure the respect of the partial ordering of tasks. Model (\ref{salbp_1} - \ref{salbp_5}) assumes all SALBP hypotheses. 

In the ALWABP,  the hypothesis of equally equipped workstations are relaxed by considering different workers in each station. \citename{miralles07advantages}~\citeyear{miralles07advantages} have proposed the ALWABP motivated by the situation found in sheltered work centers for the disabled, in which assembly lines are operated by  workers with disabilities. 
This problem has fomented an intense amount of research both on the original problem \cite{blum11solving,moreira12simple,mutlu13iterative,vila14branch,borba14heuristic} and on its variants \cite{costa09job,araujo12two,moreira13hybrid,araujo14balancing}.

In order to model the ALWABP, let $N$ be a set of tasks. The execution time of each task $i\in N$ depends on the worker it is assigned to and is assumed deterministic. Indeed, $p_{wi}$ is the executing time of task $i$ by worker $w$. Tasks $i\in N$ still must respect a partial ordering. Let $W$ represent all workers available and $S$ represent the set of workstations. Defining binary variables $y_{sw}$  (equal to 1 if, and only if, worker $w \in W$ is allocated to station $s \in S$) and  variables $x_{swi}$ (equal to 1 if, and only if, task $i\in N$ is executed at station $s\in S$ by worker $w\in W$) and  a variable $C$ representing the cycle time (execution time of the most loaded workstation), the ALWABP can be written as:

\begin{equation}
\label{alwabp_1}
\textrm{Min }   C
\end{equation}
\hspace{1cm} subjected to:
\begin{alignat}{2}
\label{alwabp_2}
\sum_{s \in S} \sum_{w \in W}  x_{swi} = 1, &&  \qquad &\forall i \in N,\\
\label{alwabp_3}
\sum_{s\in S} y_{sw} = 1, && \qquad &\forall w \in W,\\
\label{alwabp_4}
\sum_{w\in W} y_{sw} = 1, && \qquad &\forall s \in S, \\
\label{alwabp_5}
\sum_{s\in S} \sum_{w\in W }  s\cdot x_{swi} \leq  \sum_{s\in S} \sum_{w\in W} s\cdot x_{swj}, && \qquad & \forall i,j \in N, i\ipre j,\\
\label{alwabp_6}
\sum_{i\in N} \sum_{w \in W}  p_{wi}\cdot x_{swi} \leq C, && \qquad &\forall s \in S,\\
\label{alwabp_7}
x_{swi} \leq y_{sw}, && \qquad & \forall s \in S, \forall w \in W, \forall i\in N,\\
\label{alwabp_8}
y_{sw} \in \{0,1\}, && \quad &\forall s\in S, \forall w\in W,\\
\label{alwabp_9}
x_{swi} \in \{0,1\}, && \quad &\forall s\in S, \forall w\in W, \forall i\in N.
\end{alignat}

The the goal of model (\ref{alwabp_1})-(\ref{alwabp_9}) is to obtain worker and task assignments to the workstations so that the  cycle time is  minimized. Constraints (\ref{alwabp_2}) force every task to be executed, while (\ref{alwabp_3}) and (\ref{alwabp_4}) guarantee that each worker is assigned to a workstation and each workstation receives a single worker, respectively. The task execution partial ordering is respected due to constraints (\ref{alwabp_5}), while (\ref{alwabp_6}) define the line's cycle time. Finally, constraints (\ref{alwabp_7}) ensure that a task is executed by a worker in a given workstation only if the worker is assigned to the workstation.

\subsection{The bowl phenomenon and ALBPs}

Manly in non-automatized assembly lines, where people, not machines, are responsible for doing the tasks, considering execution times to be deterministic might be a strong assumption. In these cases, more precise production rates can be modelled if one considers the statistical distribution dictating the task times. We consider the case of unpaced assembly lines where the pace dictated by the cycle time is relaxed. In such unpaced lines, workstations may have to wait for the following station to finish its work in order to pass the product along the line, in this case, the station waiting is said to be blocked. On the other hand, when a station is free and has to wait for the product from the previous station, it is said to be starved. In the literature, it is common to assume that the first station can never be starved and the last one can never be blocked.

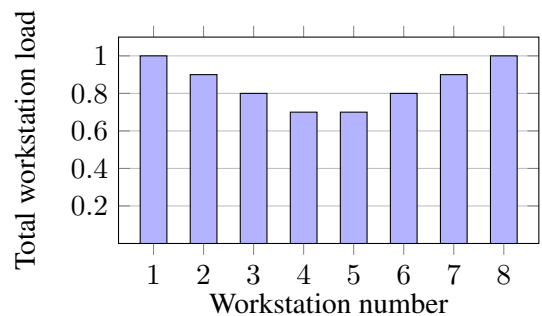
\begin{figure}[h]
\centering
\begin{tikzpicture}
\begin{axis}[ymin=0, ymax=1.1, 
	ymajorgrids,
	enlargelimits=auto,
	width=2.8in, height=1.7in,
	ylabel=Total workstation load, 
	xlabel=Workstation number,
	xtick=data,
	ytick={0.2, 0.4, 0.6, 0.8, 1.0},
	ybar
	]

\addplot[fill=blue!30!white] plot coordinates{
	(1, 1)
  	(2, 0.9)
  	(3, 0.8)
  	(4, 0.7)
  	(5, 0.7)
  	(6, 0.8)
  	(7, 0.9)
  	(8, 1.0)};
	\end{axis}
\end{tikzpicture}
\caption{example of a bowl profile. Stations in the center are less loaded than the peripheral ones.}
\label{fig:bowlShapeExample}
\end{figure}

According to \citename{Smunt1985}~\citeyear{Smunt1985}, the first studies of unpaced lines date from 1962 when simulations evaluated initial buffer configurations \cite{Baren1962}. Usually, analytical models of these kind of systems are complex and are only able to deal with no more than four workstations. Even then, they require simplifying hypotheses like assuming that the task times are exponentially distributed. More practical scenarios (with more stations and appropriated distributions) have been studied using simulation based approaches. \citename{hillier66effects} \citeyear{hillier66effects}, after some claims that the optimum allocation would be to alternate fast and slow stations \cite{Patterson1964}, showed that if the operation times (task times) were following an exponential distribution, then the assembly line throughput could be enhanced by introducing specific small load unbalances. This effect became known as the \emph{bowl phenomenon}, since the graphic of total time of each workstation along the line has the shape of a bowl, with times increasing from central station(s). Figure \ref{fig:bowlShapeExample} shows an example of such configuration in a line with eight workstations.

Results presented by \citename{hillier66effects} \citeyear{hillier66effects} considered lines with at most four stations with buffer capacity from one to four units. However, using the exponential distribution does not allow one to tell if the enhancement in the line's throughput is due to unbalancing the mean or the variability of the task times, since both parameters are the same for that distribution. The investigation of this gap, also considering exponential distribution and assuming some deterministic variables, showed the optimum strategy to be a combination of both unbalances \cite{Rao1976}. The use of \textit{Erlang} distribution also produces a bowl load profile as optimum solution and \citename{Hillier1979} \citeyear{Hillier1979} concluded that the effect of the bowl configuration became larger as the number of stations increased. By this time, the first experiments supporting the bowl phenomenon as proposed by \citename{hillier66effects} \citeyear{hillier66effects} were performed.

At one moment, conflicting results related to the production rate and the idle time of the workstations were presented. Some of them
argued that an increasing workload profile would be better than a bowl shaped workload. In this sense, \cite{Kottas1981} proposed an
approach for studying and monitoring the transient behavior of an assembly line and showed that the conflicting results were, in fact, 
underestimated transient properties read as steady-state behavior.

Another controversy happened when \citename{Smunt1985} \citeyear{Smunt1985}, using some results from \cite{Dudley1963}, argued that the assembly lines considered in previous researches, with few stations and high variability, were not of practical interest. They then presented results for more realistic scenarios that supported a degradation of the line's throughput if a bowl workload profile was used. Nevertheless, \cite{Karwan1989} and \cite{So1989}, independently argued that \cite{Smunt1985} had some flaws in their  experimental design. All the group of researchers agreed that more realistic scenarios should be explored.

\citename{Hillier1996}~\citeyear{Hillier1996} examined the robustness of the bowl phenomenon, namely, the effect of incorrect estimate the optimum bowl allocation of workload. They concluded that even errors of $50\%$ in the optimum allocation are still better than the perfectly balanced load. Moreover, a $10\%$ error in the workstations' load produces better results than the balanced case. Therefore, the authors concluded that is better to aim at a bowl profile than at a balanced one.

\citename{Hillier2006}~\citeyear{Hillier2006} optimized workload and buffers capacity simultaneously and the results showed that, if the buffers are small, their optimum allocation is balanced and the workload should follow a bowl pattern. Considering one wants to minimize total time, idle time or mean time of a product in the system, \cite{Das2010} and \cite{Das2010a}, using simulation, explored bowl, inverted bowl and crescent workload profiles. Also, in \cite{Das2012} there is a comparison between crescent and decrescent load profile. Furthermore, optimizing buffer capacity between stations and work allocation using exact and heuristic methods was studied by Hillier \citeyear{Hillier2013}.

It is important to state that none of the articles cited above considered the indivisibility of the tasks, even though this is a strong requirement in most practical situations. In other words, all authors have assumed that a task with a specific size is always available to fulfill the required workload of a station. This flaw has been already pointed out by  \citename{Tempelmeier2003} \citeyear{Tempelmeier2003} but, to the best of our knowledge, no investigation was conducted to settle this question. This is one of the main research gaps explored in the remainder of this article.

In general, the innovation presented in this paper is twofold. First, we use a simulation model to try to conclude on the existence of the bowl effect in larger lines and with  more diversified configurations. For this, we use the set of instances recently proposed by \citename{otto13salpbgen} \citeyear{otto13salpbgen} while considering the integrality of the tasks. We also propose a simple MIP methodology to generate bowl format configurations for lines with homogeneous and heterogeneous workforce, either by introducing load unbalances or by analysing the effect of unbalances in deviation. These two cases are described in the following. 

\subsubsection{Load unbalance}

\begin{figure}[h]
\centering

\pgfplotsset{}
\begin{tikzpicture}
\begin{axis}[
	title=$\alpha \mathtt{=} 1.0$  $\beta \mathtt{=} 0.9$,
	ymin=0, ymax=1.1, 
	width=2.0in, height=1.7in,
	ymajorgrids,
	enlargelimits=auto,
	ylabel=Relative workload,
	xlabel=Workstation number,
	xtick=data,
	ytick={0.2, 0.4, 0.6, 0.8, 1.0},
	bar width = 7pt,
	ybar
	]

\addplot[fill=blue!30!white] coordinates{
	(1, 1) (2, 0.9) (3, 0.81) (4, 0.9) (5, 1)};
\end{axis}

\end{tikzpicture}
\begin{tikzpicture}
\begin{axis}[
	title=$\alpha \mathtt{=} 1.0$  $\beta \mathtt{=} 0.8$,
	ymin=0, ymax=1.1, 
	width=2.0in, height=1.7in,
	ymajorgrids,
	enlargelimits=auto,
	xlabel=Workstation number,
	xtick=data,
	ytick={0.2, 0.4, 0.6, 0.8, 1.0},
	bar width = 7pt,
	ybar
	]

\addplot[fill=blue!30!white] coordinates{
	(1, 1) (2, 0.8) (3, 0.64) (4, 0.8) (5, 1)};

\end{axis}
\end{tikzpicture}
\begin{tikzpicture}
\begin{axis}[
	title=$\alpha \mathtt{=} 0.9$  $\beta \mathtt{=} 0.8$,
	ymin=0, ymax=1.1, 
	width=2.0in, height=1.7in,
	ymajorgrids,
	enlargelimits=auto,
	xlabel=Workstation number,
	xtick=data,
	ytick={0.2, 0.4, 0.6, 0.8, 1.0},
	bar width = 7pt,
	ybar
	]

\addplot[fill=blue!30!white] coordinates{	
	(1, 0.9) (2, 0.72) (3, 0.576) (4, 0.72) (5, 0.9)};	
	
\end{axis}
\end{tikzpicture}
\caption{load configurations for different values of parameters $\alpha$ and $\beta$.}
\label{fig:parametros}
\end{figure}
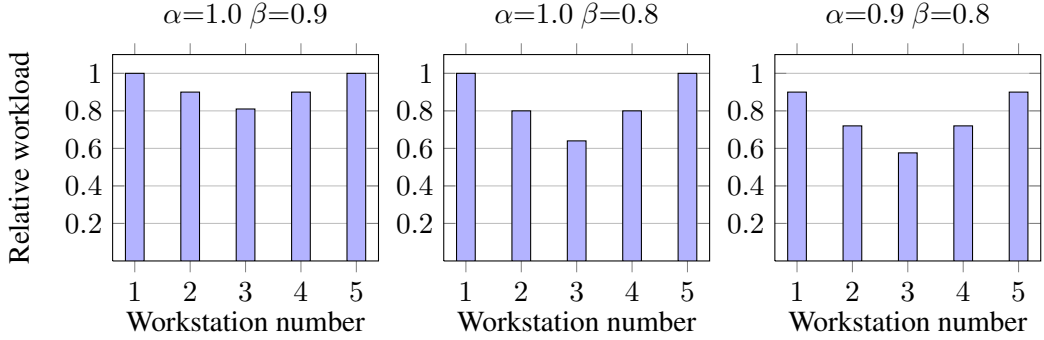

Constraints (\ref{salbp_3_cicle}) and (\ref{alwabp_6}) are appropriate for the deterministic scenario because they induce a uniform load distribution among all workstations, hence minimizing the cycle time $C$. In order to obtain bowl-like solutions, these constraints are modified in the corresponding models. For SALBP the new constraints are:

\begin{equation}
\label{salbpmodificado_1}						
\sum_{i\in N} t_{i}\cdot x_{is} \leq \alpha_s C,  \qquad  \forall s \in S,\\
\end{equation}

and for ALWABP, the cycle time constraints are replaced by:
\begin{equation}
\label{alwabpmodificado_1}						
\sum_{i\in N} \sum_{w \in W}  p_{wi}\cdot x_{swi} \leq \alpha_s C,  \qquad  \forall s \in S,\\
\end{equation}

in both cases, $\alpha_s$ are parameters whose values respect the following equations:
\begin{alignat}{2}	
\label{alwabpmodificado_2}
\alpha_1 > 0 , && \\
\label{alwabpmodificado_3}
\alpha_s  =  \alpha_{|S|-s+1},   && \qquad & s = 1 \hdots \lfloor |S|/2 \rfloor. \\
\label{alwabpmodificado_4}
\alpha_{s}  =  \beta \alpha_{s-1},   && \qquad & s = 2 \hdots \lfloor |S|/2 \rfloor - 1.   
\end{alignat}
 Constraints 
(\ref{alwabpmodificado_2}) guarantee positive loads for workstations while  constraints (\ref{alwabpmodificado_3}) induce symmetry in the  load configuration. Finally, constraints (\ref{alwabpmodificado_4}) define the bowl slope level which is controlled by parameter  $\beta$, $0 < \beta  \leq 1$. Note that if $\beta = 1$ the proposed models are equivalent to the original ones. Figure \ref{fig:parametros} shows examples of resulting profiles for some values of $\alpha$ and $\beta$.

\subsubsection{Deviation unbalance}
\label{sec:deviationUnbalance}

Other characteristic that may induce a bowl effect is related to task times deviations. Let $c_s =\frac{\sigma_s}{\mu_s}$ indicate the level of variability in a given station execution's time $s\in S$. A bowl profile is obtained by maintaining the station loads fixed and reducing the task deviation times in central stations.  The deviation was computed considering solutions without restrictions (\ref{salbpmodificado_1}) or (\ref{alwabpmodificado_1}) and by artificially adjusting task deviation times from the base case. 

The base case considers $\frac{\sigma}{\mu} = 0.1$, $s=1,...,S$, while bowl profiles are obtained with $c_s = \theta c_{s-1}$, $s= 1, ..., \lfloor |S|/2 \rfloor$ and $ c_s = c_{|S|-s+1}$, $s=1, ..., \lfloor |S|/2 \rfloor$, where $0 < \theta \leq 1$ is a parameter related to the depth of the bowl. Note that these equations are similar to the equations to generate bowl solutions for the mean case. Hence, the bowl profiles for mean and deviation have the same shapes.

\section{Simulation Model}
\label{sec:modelo}

Simulation is a technique for performance and reaction analysis of a system. It can be useful in planning stages or when the system is already in operation if careful validation on practical contexts is desired \cite{Leal2011}. Most benefits of simulation appear in contexts in which the system under study has many features that  can not be modelled properly with the use of other analytical techniques due, for example, to its complexity. 

When building a simulation model, it is important to define the inputs and outputs so the results may be analyzed. This analysis depends on the nature of the model. If the model is deterministic, for example, only one run is needed for each input one may want to analyze; on the other hand, in a stochastic model, it is possible to use a Monte Carlo method: the essence of it is to use stochastic distributions, each one representing a process in the model. Then, a statistical sampling is used to summarize the result for each input. In other words, simulation models can be used to analyze processes having a degree of uncertainty (randomness) in their variables \cite{montecarlo08}.

In this work, the simulation model built in order to evaluate assembly lines has the number of products which left the line as a state-variable, it is of continuous time, discrete state, stochastic, linear, dynamic, open and unstable \cite{Ross06Simulation}. Also, it is a next-event analysis model whose events are triggered when a station finishes its work. 

\begin{figure}[here]
\center
\includegraphics[scale=0.6]{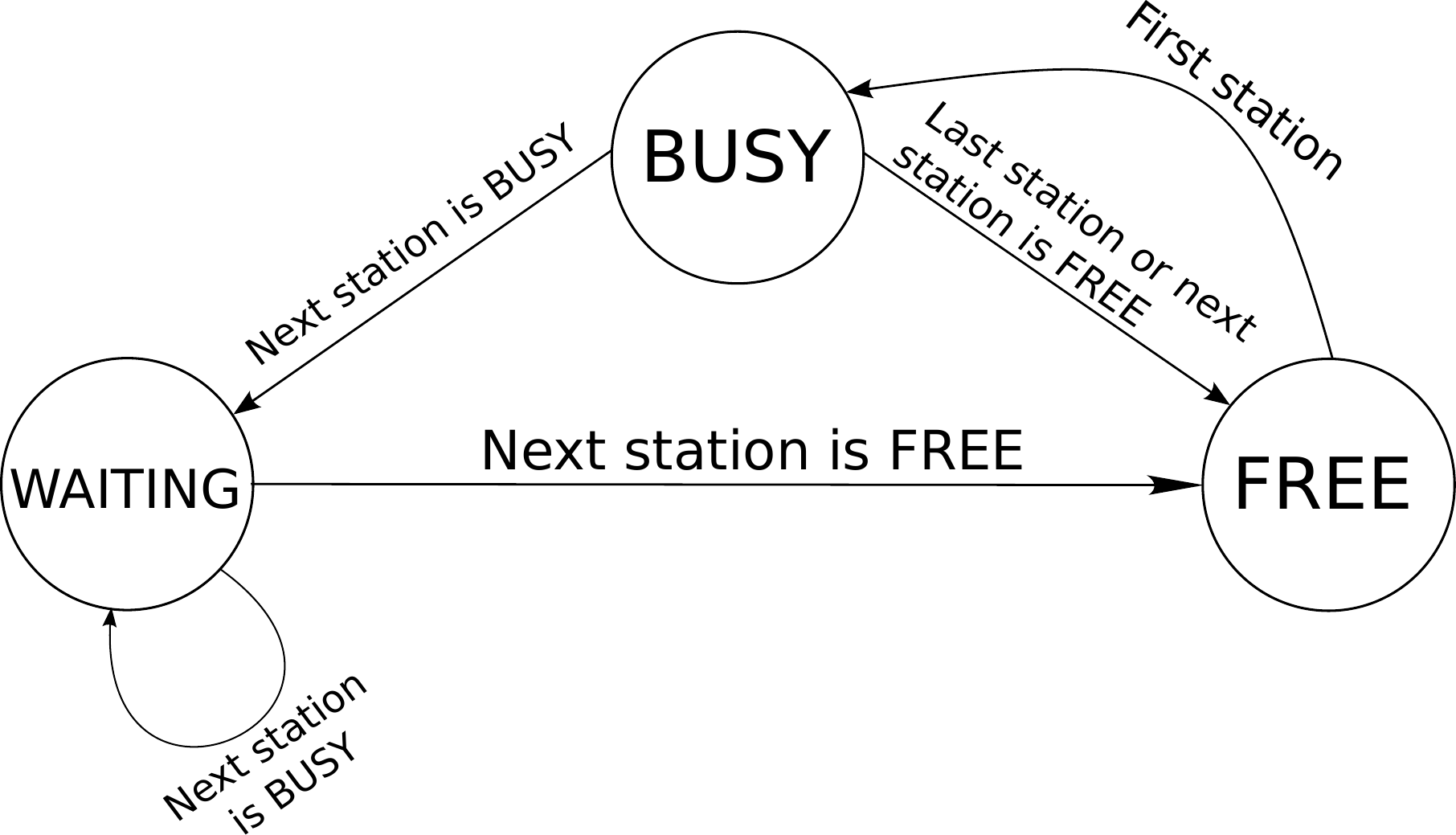}
\caption{simplified machine state for the workstations. When a workstation is in BUSY state it is performing some work, in WAITING state it is blocked and in FREE it is starved.}
\label{fig:stationStates}
\end{figure}

Actually, besides the starting event, there is just one more type of event that happens during a simulation: when all the work in a station has been done, a previously scheduled event is triggered and the product is delivered to the next station or the current station gets blocked, following the diagram in figure \ref{fig:stationStates}. There are three possible workstation states. Any transition can only happen when the referred event is triggered. Namely, whenever the event happens, if the current station is in BUSY state it checks the following station, if it is in FREE state, the finished piece of product is delivered to it and the state of the current station changes to FREE, which means it can receive another piece to work on, if it is the first station it returns to BUSY state (there is always work available for the first station); otherwise, if the next station is in BUSY state when the current finishes its work, the current station changes its state to WAITING, blocking itself until the next station is FREE. If the current BUSY station is the last station of the line it changes its state to FREE, since there is no need to check for the next station's state. There is also one action, not shown in figure \ref{fig:stationStates},  that takes place when the transition from BUSY to FREE happens: the current station changes the state of the next station to BUSY, translating the fact that a new piece of product arrived in the following station. Hence, the event triggered when a workstation finishes its work controls the flow of products in the line. 

\begin{figure}[here]
\center
\includegraphics[scale=0.6]{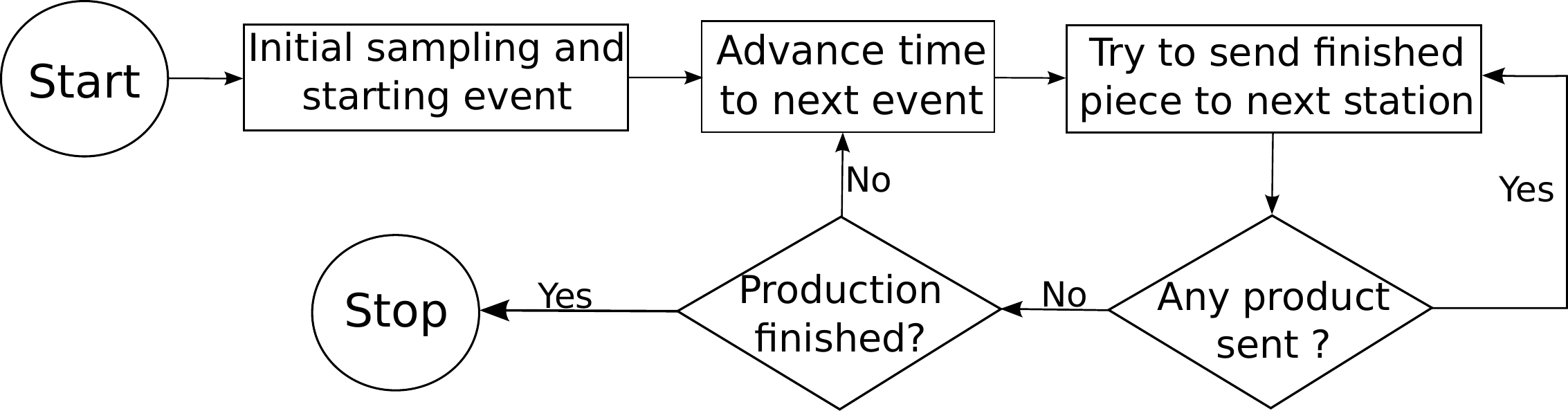}
\caption{flowchart of a simulation. The simulation advances its time to the next scheduled event and process the corresponding events. If none of them can be processed it advances time again until at least one station is not blocked. Simulation stops when a specified number of completed products is reached.}
\label{fig:simulationModel}
\end{figure}

The scheduling of an event occurs when the station changes its state from FREE to BUSY. In the model, this transition starts a sampling procedure of the station tasks execution times and the event associated with the completion of the work in that station is scheduled to the proper time.  The handling of scheduled events is performed according to the diagram shown in figure \ref{fig:simulationModel} which has been adapted from \cite{robinson2004}.

When a simulation starts, every station state is set to FREE, but the first one whose state is set to BUSY, hence an event corresponding to finishing the first job in the first station is scheduled. Then, time is advanced to the next event which is processed according to diagram in figure \ref{fig:stationStates}. If no station can pass the product to the next station and the desired number of products is not yet reached, the simulation advances to the next scheduled event until a workstation becomes able to deliver its work to the following one. The simulation stops when the specified production is reached. 

The allocation of tasks to workstations and, in ALWABP case, workers to workstations is given by the respective deterministic mathematical model whose objective is to minimize the assembly line's cycle time. Since, in the simulation model of this work, every task is supposed to be performed in a time following a Normal distribution, two parameters are needed in order to fully specify the task times to the simulation model: mean $\mu$ and standard deviation $\sigma$. The mean is assumed to be the deterministic time of the task and the standard deviation is $\sigma = \frac{\mu}{10}$ for load unbalance simulations and is given according to subsection \ref{sec:deviationUnbalance} for deviation unbalance simulations. 

In summary, we propose a  stochastic simulation model using next-event analysis and deterministic mathematical programming models  to study the bowl phenomenon considering integer tasks and instances proposed by \cite{otto13salpbgen} and \cite{chaves09hybridb}, respectively for SALBP and ALWABP. The results obtained in the experiments are detailed in the next section.

\section{Computational results}
\label{sec:resultados}

The computational results were obtained with a series of stochastic simulations using a Monte Carlo method. Each instance was run 300 times in order to obtain statistically relevant results, and each simulation run until a production of 150 items was reached. The time to accomplish this production is, therefore, the result to be analyzed. For each run the results obtained for the first 50 items were discarded in order to account for a possible transient state of the line. This warm-up time was obtained using Welch's method as described in \cite{Mahajan2004}.  All the task times were assumed to follow a Normal distribution, $N(\mu, \sigma^2)$, and separated experiments were performed to evaluate the  mean and standard deviation unbalances on total workstation load. In order to analyze mean unbalances, experiments with $\alpha_1 = 1$ and $\beta \in B = \{1.00, 0.99, 0.98, 0.97, 0.96, 0.95, 0.94\}$ were performed (as defined in constraints (\ref{alwabpmodificado_2}) and (\ref{alwabpmodificado_4})); for deviation unbalances the experiments considered $\theta  \in \{\theta_i = {\beta_i} | \forall \beta_i \in B\}$.

Both SALBP and ALWABP cases were considered independently. The benchmark problems for the SALBP simulations are a subset of the   dataset proposed by \cite{otto13salpbgen}. We use the instances for which the known upper bound on the number of stations is seven or less. The benchmark for the ALWABP simulations are the instances of type Heskia and Roszieg from \cite{chaves09hybridb}. Both, SALBP and ALWABP instances, were solved to optimality. The results' analysis, performed in the following, focus on the SALBP problems since their instances are more representative than their ALWABP counterparts. Tables \ref{tab:SALBPmean_ws} through \ref{tab:SALBPdev_ws_gs_td} summarize results for the SALBP simulations. Each table presents the total number of instances in the corresponding categories, the number of instances in which the bowl phenomenon could be observed (i.e., the number of instances in which the bowl configuration was statistically superior to their balanced counterparts) and their mean $\beta$ (or $\theta$) for mean (or standard deviation) unbalances. The categories are related to different precedence graph structures and how the task times are distributed: graphs may have \textit{bottleneck} tasks and \textit{chain} tasks and, according to the relevance of one or other type of tasks, \citename{otto13salpbgen} \citeyear{otto13salpbgen} classified the graph structure as \textit{BN} (for \textit{bottleneck}), \textit{CH} (for \textit{chain}) or \textit{MIXED} when there is no prominent characteristic. Task time's distributions may be bimodal or with peak at small tasks (bottom peak) in the considered instances. Also a statistical t-test with $0.05$ threshold in the p-value was used to ensure statistical difference from the best bowl solution to the balanced one.

\begin{table}[!ht]
\centering
\caption{mean unbalance results for SALBP instances grouped by number of workstations.}

\begin{tabular}{|c|c|c|c|c|c|c|c|c|c|}
\hline
& \multicolumn{5}{c|}{\textbf{Workstations}}\\
\hline
& \textbf{3} & \textbf{4} & \textbf{5} & \textbf{6} & \textbf{7}\\
\hline
\textbf{Instances} & 163 & 38 & 89 & 50 & 6\\ 
\hline
\textbf{Bowl} & 0.5399 & 0.8421 &0.7416 & 0.7800 & 0.8333 \\
\hline
\textbf{$\beta_{mean}$} & 0.9781	& 0.9781 	& 0.9788 	& 0.9777 & 	0.9900\\
\hline

\end{tabular}
\label{tab:SALBPmean_ws}
\end{table}

Tables \ref{tab:SALBPmean_ws}, \ref{tab:SALBPmean_gs_td} and \ref{tab:SALBPmean_ws_gs_td} contain results related to the mean unbalance simulations. Table \ref{tab:SALBPmean_ws} presents the results grouped by number of stations. Most of the instances ($163$) considered have three workstations in the optimal solution and in $53.99\%$ the best bowl solution was statistically different from the balanced solution leading to $\beta_{mean} = 0.9781$. The deepest mean bowl was also found in instances with four workstations which was, as well, the case that present greater relative number of observable bowl phenomenons ($84.21\%$).

\begin{table}[!ht]
\centering
\caption{mean unbalance results for SALBP instances grouped by graph structure and task time distribution.}
\begin{tabular}{|c|c|c|c|c|c|}
\hline
& \multicolumn{3}{c|}{\textbf{Graph structure}} & \multicolumn{2}{c|}{\textbf{Task time distribution}}\\
\hline
& \textbf{BN} & \textbf{CH} & \textbf{MIXED} & \textbf{Bimodal} & \textbf{Bottom peak}\\
\hline
\textbf{Instances}  & 98 & 99 & 149 & 174 & 172\\ 
\hline
\textbf{Bowl} & 0.6531	& 0.6970 & 	0.6510 &	0.7644&	0.5640\\
\hline
\textbf{$\beta_{mean}$}	&  0.9820 & 0.9781	& 0.9764 & 0.9788 & 0.9780\\
\hline

\end{tabular}
\label{tab:SALBPmean_gs_td}
\end{table}

Table \ref{tab:SALBPmean_gs_td} summarizes the results grouped by graph structures and task time distributions. In the graph structure group, the deepest mean bowl was found in \textit{MIXED} and the \textit{CH} class presented the greater relative number of bowl-shaped solutions performing better than balanced one. The bottom peak task time distribution presented the deepest bowl and bimodal instances shown a greater number of solutions for which the bowl configuration was more productive. Also, it can be seen that \textit{CH} graphs and bimodal task distribution seem to favor the occurrence of the bowl phenomenon.

\begin{table}[!ht]
\centering
\caption{mean unbalance results for SALBP instances combining number of workstation with graph structure or task time distribution.}
\begin{tabular}{|c|c|c|c|c|c|c|}
\hline
& & \multicolumn{3}{c|}{\textbf{Graph structure}} & \multicolumn{2}{c|}{\textbf{Task time distribution}}\\
\hline
\textbf{Workstations} & & \textbf{BN} & \textbf{CH} & \textbf{MIXED} & \textbf{Bimodal} & \textbf{Bottom peak}\\
\hline
\multirow{3}{*}{\textbf{3}} & \textbf{Instances}  &  47 & 46 & 70 & 0  & 163 \\
\hhline{~------}		       
&	 \textbf{Bowl}	& 0.5319	 & 0.5435 & 0.5429	 & - & 	0.5399 \\
\hhline{~------}
&\textbf{$\beta_{mean}$} & 0.9832	& 0.9772	 & 0.9753	& - &	0.9781 \\
\hline
\multirow{3}{*}{\textbf{4}} & \textbf{Instances}  &  11 & 12 & 15 & 29  & 9 \\
\hhline{~------}		       
&	 \textbf{Bowl}	& 0.6364 &	1.0000 &	0.8667	& 0.7931 	& 1.0000 \\
\hhline{~------}
&\textbf{$\beta_{mean}$}& 0.9800	&0.9792&	0.9762&	0.9783	&0.9778 \\
\hline
\multirow{3}{*}{\textbf{5}} & \textbf{Instances}  &  28 & 24 & 37 & 89  & 0 \\
\hhline{~------}		       
&	 \textbf{Bowl}	& 0.7857 &	0.7083&	0.7297	&0.7416	&- \\
\hhline{~------}
&\textbf{$\beta_{mean}$} & 0.9818 &	0.9770 & 	0.9774	&0.9788	&- \\
\hline
\multirow{3}{*}{\textbf{6}} & \textbf{Instances}  &  10 & 15 & 25 & 50  & 0 \\
\hhline{~------}		       
&	 \textbf{Bowl}	& 0.8000	 & 0.8667 &	0.7200	& 0.7800 	& - \\
\hhline{~------}
&\textbf{$\beta_{mean}$} & 0.9788	& 0.9785	& 0.9767	& 0.9777	& - \\
\hline
\multirow{3}{*}{\textbf{7}} & \textbf{Instances}  &  2 & 2 & 2 & 6  & 0 \\
\hhline{~------}		       
&	 \textbf{Bowl} &  1.0000	& 1.0000	& 0.5000 & 0.8333 	& - \\
\hhline{~------}
&\textbf{$\beta_{mean}$} & 0.9900 & 0.9900	 & 0.9900 & 0.9900 & - \\
\hline
\end{tabular}
\label{tab:SALBPmean_ws_gs_td}
\end{table}

Table \ref{tab:SALBPmean_ws_gs_td} presents the results considering both number of stations and graph structure and number of stations and task time distribution. In some of these scenarios there was no available instance: three workstations with bimodal task time distribution and five workstations with bottom peak time distribution, for example. Therefore, there are table cells marked with "-" where the corresponding result is not applicable. Again, for three and four workstation, it is possible to see \textit{CH} graphs favouring the bowl phenomenon, nevertheless, this is not true for five workstations. It is hard to draw conclusions for instances with seven workstations since there are fewer instances in this scenario. Likewise, for task time distribution case the instances are not well suited for comparing bimodal and bottom peak distributions.

\begin{table}[!ht]
\centering
\caption{deviation unbalance results for SALBP instances grouped by number of workstations.}

\begin{tabular}{|c|c|c|c|c|c|c|c|c|c|}
\hline
& \multicolumn{5}{c|}{\textbf{Workstations}}\\
\hline
& \textbf{3} & \textbf{4} & \textbf{5} & \textbf{6} & \textbf{7}\\
\hline
\textbf{Instances} & 163 & 38 & 89 & 50 & 6\\ 
\hline
\textbf{Bowl} & 0.5337	& 0.9737	& 1.0000	& 1.0000	& 1.0000 \\
\hline
\textbf{$\theta_{mean}$} & 0.9480	& 0.9462 &	0.9402 &	 0.9410	& 0.9400\\
\hline

\end{tabular}
\label{tab:SALBPdev_ws}
\end{table}

Tables \ref{tab:SALBPdev_ws}, \ref{tab:SALBPdev_gs_td} and \ref{tab:SALBPdev_ws_gs_td} present results regarding bowl profiles in $c_s = \frac{\sigma_s}{\mu_s}$. Table \ref{tab:SALBPdev_ws} shows results grouped by number of workstations. Instances with five, six and seven workstations achieved $100\%$ of \emph{bowl phenomenon occurrence} (i.e., in all cases the lines with the obtained bowl configurations were more productive than their balanced counterparts). Also, the mean depth of the bowl for seven workstation reached the maximum value among all the experiments conducted. 

\begin{table}[!ht]
\centering
\caption{deviation unbalance results for SALBP instances grouped by graph structure and task time distribution.}
\begin{tabular}{|c|c|c|c|c|c|}
\hline
& \multicolumn{3}{c|}{\textbf{Graph structure}} & \multicolumn{2}{c|}{\textbf{Task time distribution}}\\
\hline
& \textbf{BN} & \textbf{CH} & \textbf{MIXED} & \textbf{Bimodal} & \textbf{Bottom peak}\\
\hline
\textbf{Instances}  & 98 & 99 & 149 & 174 & 175\\ 
\hline
\textbf{Bowl}& 0.7653 &	0.8080 &	0.7651& 	0.9943	& 0.5581\\
\hline
\textbf{$\theta_{mean}$} & 0.9437 & 0.9446 &	0.9445 &	0.9424 &	0.9478\\
\hline

\end{tabular}
\label{tab:SALBPdev_gs_td}
\end{table}

Table \ref{tab:SALBPdev_gs_td} presents results grouped by graph structures and task time distributions. In the graph structure group, the deepest mean bowl was found in the \textit{bottleneck} class and the class that presented the greater relative number of bowl-shaped solution performing better than balanced solutions was the \textit{chain} one. The bimodal task time distribution presented deepest and greater number of better bowl solutions. Four out of five classes presented more than $75\%$ of occurrence of the bowl phenomenon, reaching $99.43\%$ in the bimodal class.

\begin{table}[!ht]
\centering
\caption{deviation unbalance results for SALBP instances combining number of workstation with graph structure or task time distribution}
\begin{tabular}{|c|c|c|c|c|c|c|}
\hline
& & \multicolumn{3}{c|}{\textbf{Graph structure}} & \multicolumn{2}{c|}{\textbf{Task time distribution}}\\
\hline
\textbf{Workstations} & & \textbf{BN} & \textbf{CH} & \textbf{MIXED} & \textbf{Bimodal} & \textbf{Bottom peak}\\
\hline
\multirow{3}{*}{\textbf{3}} & \textbf{Instances}  &  47 & 46 & 70 & 0  & 163 \\
\hhline{~------}		       
&	 \textbf{Bowl}  & 0.5319  &	0.5869  &	0.5000 &	-	 &0.5337 \\
\hhline{~------}
&\textbf{$\theta_{mean}$} & 0.9472	& 0.9478 	& 0.9489 &	-	&0.9480 \\
\hline
\multirow{3}{*}{\textbf{4}} & \textbf{Instances}  &  11 & 12 & 15 & 29  & 9 \\
\hhline{~------}		       
&	 \textbf{Bowl} & 0.9090	& 1.0000	& 1.0000	& 0.9655 & 1.0000 \\
\hhline{~------}
&\textbf{$\theta_{mean}$} & 0.9430	& 0.9492 & 	0.9460 & 	0.9464	& 0.9456 \\
\hline
\multirow{3}{*}{\textbf{5}} & \textbf{Instances}  &  28 & 24 & 37 & 89  & 0 \\
\hhline{~------}		       
&	 \textbf{Bowl}		&  1.0000  & 1.0000 & 1.0000 &  1.0000 &  - \\
\hhline{~------}
&\textbf{$\theta_{mean}$} & 0.9421	& 0.9413	 & 0.9424	& 0.9420	 & - \\
\hline
\multirow{3}{*}{\textbf{6}} & \textbf{Instances}  &  10 & 15 & 25 & 50  & 0 \\
\hhline{~------}		       
&	 \textbf{Bowl}		&  1.0000 & 1.0000 & 1.0000 & 1.0000 &  - \\
\hhline{~------}
&\textbf{$\theta_{mean}$} & 0.9410	& 0.9413	 & 0.9408 &	0.9410 & - \\
\hline
\multirow{3}{*}{\textbf{7}} & \textbf{Instances}  &  2 & 2 & 2 & 6  & 0 \\
\hhline{~------}		       
&	 \textbf{Bowl}		&  1.0000 & 1.0000 & 1.0000 & 1.0000 &  - \\
\hhline{~------}
&\textbf{$\theta_{mean}$} & 0.9400 & 0.9400 & 0.9400  & 0.9400 & - \\
\hline

\end{tabular}
\label{tab:SALBPdev_ws_gs_td}
\end{table}

Table \ref{tab:SALBPdev_ws_gs_td} summarizes results combining number of workstations and graph structure and number of workstations and task time distribution for deviation analysis. For instances with three workstations the graph structure \textit{chain} seems to favor the effect. For instances with more than three stations the phenomenon could be observed in most of the problems and the maximum bowl depth considered in the experiments was in a number of scenarios. 

\begin{table}[!ht]
\centering
\caption{mean and deviation unbalances for ALWABP instances.}
\begin{tabular}{|c|c|c|c|c|}
\hline
& \multicolumn{2}{c|}{\textbf{Low}} & \multicolumn{2}{c|}{\textbf{High}}\\
\hhline{~----}
& \textbf{mean} & \textbf{deviation} & \textbf{mean} & \textbf{deviation}\\
\hline
\textbf{Bowl} & 0.6125 & 0.8625 & 0.7500 & 0.8000 \\
\hline
$\beta_{mean}$ or $\theta_{mean}$ & 0.9684 & 0.9423 & 0.9670 & 0.9430	\\
\hline
\end{tabular}
\label{tab:ALWABPresults}
\end{table}

Regarding the ALWABP, the results were analyzed in a much simpler way since the available instances are not as representative as those for the SALBP. Table \ref{tab:ALWABPresults} shows the results. They are grouped by variance among workers: low variance means that  
task execution time varies less among workers than in the high variance scenario \cite{chaves09hybridb}. The table presents the fraction of instances where the bowl phenomenon could be observed and $\beta_{mean}$ for mean unbalance or $\theta_{mean}$ for deviation unbalance. For each 
class (low and high) $80$ instances were analyzed. It is possible to see that $\beta_{mean}$ had smaller values here than in the  SALBP case. Furthermore, in the deviation case the bowl phenomenon was observed in fewer instances if compared to SALBP. Nevertheless, the bowl phenomenon could be observed in more than $50\%$ of the cases for every class (and case) considered.

As a final note, we remark that the SALBP and ALWABP results show that the bowl phenomenon regarding variability occurred in more instances than the effect regarding mean. This may have happened because there is no integrality constraint in the standard deviation whereas, for the mean case, the proposed MIP formulation considers the integrality of the tasks.

\section{Conclusions}
\label{sec:conclusoes}

This study introduced mixed-integer mathematical models to obtain slightly unbalanced load profiles for the simple assembly line balancing problem and the assembly line worker assignment balancing problem. These profiles were used in a simulation model to verify that the efficiency of the lines could be improved if they were slightly unbalanced. This was done using  recently proposed benchmark problems, which are more representative than those used in previous bowl phenomenon studies, favoring the generalization of the conclusions obtained. These conclusions  indicate that assembly lines can indeed benefit from bowl shaped configurations even if the more realistic scenario with indivisible tasks is considered.

\section{Acknowledgments}

This work was developed with financial support from CNPq, CAPES/DGU (258-12) and FAPESP.

\bibliography{castellucci2012eng}

\end{document}